\documentclass[12pt]{article}
\usepackage{amsmath,amsthm,amsfonts,amssymb,times}

\theoremstyle{remark}
\newtheorem{remark}{Remark}
\theoremstyle{plain}

\newtheorem{lem}{Lemma}[section]
\newtheorem{thm}{Theorem}

\numberwithin{equation}{section}

\newcommand{\Z}{{\mathbb Z}}

\newcommand{\CC}{{\mathcal C}}

\newcommand{\lcm}{{\rm lcm}}
\newcommand{\be}{\begin{equation}}
\newcommand{\ee}{\end{equation}}
\newcommand{\benn}{\begin{equation*}}
\newcommand{\eenn}{\end{equation*}}

\newcommand{\lam}{\ensuremath{\lambda}}

\renewcommand{\a}{\ensuremath{\alpha}}
\renewcommand{\b}{\ensuremath{\beta}}

\newcommand{\del}{\ensuremath{\delta}}
\newcommand{\eps}{\ensuremath{\varepsilon}}
\renewcommand{\(}{\left(}
\renewcommand{\)}{\right)}

\newcommand{\pfrac}[2]{\left(\frac{#1}{#2}\right)}

\newcommand{\tC}{\widetilde{C}}

\allowdisplaybreaks

\begin{document}

\title{\vskip -40pt
\LARGE \textbf{\textsc{Sieving by large integers \\
and covering systems of congruences\\ \ }}}
\author{Michael Filaseta
\\ Department of Mathematics \\
University of South Carolina \\ Columbia, SC  29208\\
filaseta@math.sc.edu
\and
Kevin Ford
\\  Department of Mathematics \\
University of Illinois  \\
Urbana, IL 61801\\
ford@math.uiuc.edu
\and \\
Sergei Konyagin
\\ Department of Mathematics \\
Moscow State University \\ Moscow, Russia  119992\\
konyagin@ok.ru
\and \\
Carl Pomerance
\\ Department of Mathematics \\
Dartmouth College \\ Hanover, NH  03755-3551\\
carl.pomerance@dartmouth.edu
\and \\
Gang Yu
\\ Department of Mathematics \\
University of South Carolina \\ Columbia, SC  29208\\
yu@math.sc.edu
}
\date{\today}

\maketitle

{\renewcommand{\thefootnote}{}
\vskip -25pt
\footnote{
\vskip 0pt
The first  author was supported by NSF grant DMS-0207302 and
NSA grant H98230-05-1-0038.
The second author was supported by NSF grant DMS-0301083.
Much of the research for this paper was accomplished while the third
author was visiting the University of South Carolina, Columbia in
January 2004 (supported by NSF grant DMS-0200187)
and the University of Illinois at Urbana-Champaign in
February 2004 (supported by NSF grant DMS-0301083).
The fourth author was supported by NSF grant DMS-0401422.}
}

\newpage
\begin{abstract}
An old question of Erd\H os asks if there exists, for each number $N$, a
finite set $S$ of integers greater than $N$ and residue classes 
$r(n)~({\rm mod}~n)$
for $n\in S$ whose union is $\Z$.  We prove that if $\sum_{n\in S}1/n$ is
bounded for such a covering of the integers, then the least member of 
$S$ is also bounded,
thus confirming a conjecture of Erd\H os and Selfridge.  We also 
prove a conjecture
of Erd\H os and Graham, that, for each fixed number $K>1$, the complement in
$\Z$ of any union of residue classes $r(n)~({\rm mod}~n)$, for
distinct $n\in(N,KN]$, has density
at least $d_K$
for $N$ sufficiently large.
Here $d_K$ is a positive number depending only on $K$.
Either of these new results implies another conjecture of Erd\H os and Graham,
that if $S$ is a finite set of moduli greater than $N$, with a choice for
residue classes $r(n)~({\rm mod}~n)$ for $n\in S$ which covers $\Z$, then the
largest member of $S$ cannot be $O(N)$.  We further obtain stronger forms
of these results and establish other information, including an improvement
of a related theorem of Haight.

\end{abstract}

\section{Introduction}\label{intro}
Notice that every integer $n$ satisfies at least one of the congruences
\[
n\equiv0~({\rm mod}~2),~n\equiv0~({\rm mod}~3),~n\equiv1~({\rm 
mod}~4),~n\equiv1~({\rm mod}~6),~n\equiv11~({\rm mod}~12).
\]
A finite set of congruences, where each integer satisfies at 
least one them, is called a \textit{covering system}.
A famous problem of Erd\H os from 1950 \cite{E2} is to determine whether
for every $N$, there is a covering system
with distinct moduli greater than $N$.
In other words, can the minimum modulus in a covering system with
distinct moduli be arbitrarily large?   In regards to this
problem, Erd\H os writes in \cite{E3}, ``This is
perhaps my favourite problem."

It is easy to see that in a covering system, the reciprocal sum of 
the moduli is at least 1.  Examples with distinct moduli are known
with least modulus 2, 3, and 4, where this reciprocal sum can be
arbitrarily close to 1, see \cite{G}, \S F13.  Erd\H os and Selfridge
\cite{E2pt5} conjectured that this fails
for all large enough choices of the least modulus.  In fact, they
made the following much stronger conjecture.

\medskip
\noindent
{\bf Conjecture 1.}
{\em For any number $B$, there is a number $N_B$, such that in a 
covering system with
distinct moduli greater than $N_B$, the sum of reciprocals of these moduli is
greater than $B$.}
\medskip

\noindent
A version of Conjecture 1 also appears in \cite{EG}.  

Whether or not one can cover all of $\Z$, it is interesting to consider how
much of $\Z$ one can cover with residue classes $r(n)~({\rm mod}~n)$, where
the moduli $n$ come from an interval $(N,KN]$ and are distinct.  In 
this regard,
Erd\H os and Graham~\cite{EG} have formulated the following conjecture.

\medskip
\noindent
{\bf Conjecture 2.}
{\em For each number $K>1$ there is a positive number $d_K$ such that 
if $N$ is sufficiently
large, depending on $K$, and we choose arbitrary integers $r(n)$ for 
each $n\in (N,KN]$, then
the complement in $\Z$ of the union of the residue classes $r(n)~({\rm mod}~n)$
has density at least $d_K$.}
\medskip

\noindent
In \cite{E3}, Erd\H os writes with respect to establishing such a lower bound
$d_K$ for the density, ``I am not sure at all if this is possible and
I give \$100 for an answer."

A corollary of either Conjecture 1 or Conjecture 2
is the following conjecture also raised by Erd\H os and Graham in \cite{EG}.

\medskip
\noindent
{\bf Conjecture 3.}
{\em For any number $K > 1$ and $N$ sufficiently large, depending on $K$, there
is no covering system using distinct moduli from the interval $(N, K N]$.}
\medskip

In this paper we prove strong forms of Conjectures 1, 2, and 3.

Despite the age and fame of the minimum modulus problem, there are still many
more questions than answers.  We mention a few results.  
Following earlier work of Churchhouse,
Krukenberg, Choi, and Morikawa, Gibson~\cite{G1} has recently constructed
a covering system with minimum modulus 25, which stands as the
largest known least modulus for a covering system with distinct moduli.
As has been mentioned, if $r_i~({\rm mod}~n_i)$ for $i=1,2,\dots,l$
is a covering system, then $\sum1/n_i\ge 1$.  Assuming that the moduli
$n_i$ are distinct and larger than 1, it is possible to show that 
equality cannot occur, that is,
$\sum1/n_i>1$.  The following proof (of M.~Newman) is a gem.  Suppose 
that $\sum1/n_i=1$.  If the system
then covers, a density argument shows that there cannot be any overlap between
the residue classes, that is, we have an {\em exact covering system}.
We suppose, as we may, that $n_1 < n_2 < \cdots < n_l$ and
each $r_i\in[0,n_i-1]$.  Then
\[
\frac1{1-z}=1+z+z^2+\dots=\sum_{i=1}^l\(z^{r_i}+z^{r_i+n_i}+z^{r_i+2n_i}+\dots\)
=\sum_{i=1}^l\frac{z^{r_i}}{1-z^{n_i}}.
\]
The right side of this equation has poles at the primitive $n_l$-th roots of 1,
which is not true of the left side.  Thus, there
cannot be an exact covering system with distinct moduli greater than 1
(in fact, the largest modulus must be repeated).

Say an integer $H$ is ``covering" if there is a covering system with distinct
moduli with each modulus a divisor of $H$ exceeding 1.  For example, 12 is
covering, as one can see from our opening example.  From the above 
result, if $H$
is covering, then $\sigma(H)/H>2$, where $\sigma$ is the sum-of-divisors
function.  Benkoski and Erd\H os~\cite{BE} wondered if $\sigma(H)/H$ 
is large enough,
would this condition suffice for $H$ to be covering.  In \cite{H}, Haight
showed that this is not the case.  We obtain a strengthening of this result,
and by a shorter proof.

If $n_1,n_2,\dots,n_l$ are positive integers and
$C=\{(n_i,r_i):i=1,2,\dots,l\}$ is a set of ordered pairs,
let $\delta=\delta(C)$
be the density of the integers that are {\it not} in the union of 
the residue classes $r_i\pmod{n_i}$.
If $n_1,n_2,\dots,n_l$ are pairwise coprime, there is no mystery
about $\delta$.  Indeed, 
the Chinese remainder theorem implies that 
for any choice of residues $r_1,r_2,\dots,r_l$, 
\[
\delta=\prod_{i=1}^l(1-1/n_i),
\]
which is necessarily positive if each $n_i>1$.

One central idea in this paper is to determine how to
estimate $\delta $ when the moduli are not necessarily pairwise coprime.
We note that for any $n_1,n_2,\dots,n_l$, there is a choice for 
$r_1,r_2,\dots,r_l$ such that 
\be
\label{deltaineq}
\delta\le\prod_{i=1}^l(1-1/n_i).
\ee
Indeed, this is obvious if $l=1$.  Assume it is true for $l$, and say we have 
chosen residues $r_1,r_2,\dots,r_l$ such that the residual set $R$
has density $\delta$ satisfying \eqref{deltaineq}.
The residue classes modulo $n_{l+1}$ partition any
subset of $\Z$, and in particular partition $R$, so that at least one
of these residue classes, when intersected with $R$,
has density at least $\delta/n_{l+1}$.
Removing such a residue class, the residual set for the $l+1$ congruences thus
has density at most
\[
\delta-\delta/n_{l+1}=(1-1/n_{l+1})\delta\le\prod_{i=1}^{l+1}(1-1/n_i).
\]
Thus, the assertion follows.

Note that
\[
\prod_{N<n\le KN}(1-1/n)
=\lfloor N\rfloor/\lfloor KN\rfloor
\to1/K\hbox{ as }N\to\infty,
\]
so that $d_K$ in Conjecture 2 must be at most $1/K$.  
We show in Section~\ref{lowerbounds} that
any number $d<1/K$ is a valid choice for $d_K$.

A key lemma in our paper allows us to almost reverse the inequality
\eqref{deltaineq} for $\delta$.  Namely we show that for any choice
of residues $r_1,r_2,\dots,r_l$, 
\begin{equation}
\label{keylemma}
\delta\ge\prod_{i=1}^l(1-1/n_i)
-\sum_{\substack{i<j\\ \gcd(n_i,n_j)>1}}\frac1{n_in_j}.
\end{equation}
We then maneuver to show that under certain conditions the product
is larger than the sum, so that no choice of residue classes $r_i$
allows a covering.  As kindly pointed out to us by the referee,
the inequality \eqref{keylemma} bears a resemblance to the
Lov\'asz Local Lemma, but seems to be independent of it.  We
shall discuss this connection more in the next section. 

If $S$ is a finite set of positive integers, 
let $\delta^-(S)$ be the minimum value of $\delta(C)$
where $C$ runs over all choices of $\{(n,r(n)):n\in S\}$.
That is, we are given the moduli $n\in S$, and we choose the residue
classes $r(n)\pmod{n}$ so as to cover as much as possible from $\Z$;
then $\delta^-(S)$ is the density of the integers not covered.
Further, let
\[
\alpha(S)=\prod_{n\in S}(1-1/n),
\]
so that \eqref{deltaineq} implies we have $\delta^-(S)\le\alpha(S)$.
With this notation we now state our principal results.

\medskip

\noindent{\bf Theorem A.}
{\em
Let $0<c<1/3$ and let $N$ be sufficiently large (depending on $c$).
If $S$ is a finite set of integers $n>N$ such that
\[
\sum_{n\in S}\frac1n\le c \frac{\log N\log\log\log N}{\log\log N},
\]
then $\delta^-(S)>0$.}
\medskip

\noindent{\bf Theorem B.}
{\em For any numbers $c$ with $0<c<1/2$, $N\ge20$, and $K$ with
$$
1<K\le\exp(c\log N\log\log\log N/\log\log N),
$$
if $S$ is a set of integers contained in $(N,KN]$, then
\[
\delta^-(S)=(1+o(1))\alpha(S)
\]
as $N\to\infty$, where the function ``$o(1)$" depends only on the 
choice of $c$.}
\medskip

\noindent
Theorems A and B are proved in Section \ref{lowerbounds}.
Using Lemma \ref{combinedlemma} below, we can make the $o(1)$ term
in Theorem B explicit in terms of $N$ and $c$.  Both Theorems~A and B,
as well as several other of our results,
are proved in a more general context of multisets $S$ or, equivalently,
where multiple residue classes are allowed for each modulus.
Note that Theorems~A and B prove Conjectures~1 and 2, respectively, and so
Conjecture~3 as well.

In the context of Theorem B, if we relax the upper bound on the largest
modulus, we are able to construct examples of sets of integers $S$
with least member arbitrarily large and where $\delta^-(S)$ is much
smaller than $\alpha(S)$.  Proved in Section~\ref{near}, this result 
might be interpreted as lending weight towards the existence
of covering systems with the least modulus being arbitrarily large.

Similar to the definition of $\delta^-(S)$, let $\delta^+(S)$ be the
{\em largest} possible density
for a residual set with $S$ being a set of (distinct) moduli.  
It was shown by Rogers, see \cite{HRo}, pp.~242--244, 
that for any finite set of positive integers $S$, the density
$\delta^+(S)$ is attained when we choose the residue class 0~(mod~$n$)
for each $n\in S$.  
That is, $\delta^+(S)$ is the density of integers not divisible
by any member of $S$.  There is an extensive literature on estimating
$\delta^+(S)$ when $S$ consists of all integers in an interval
(see e.g.\ \cite{F} and Chapter 2 of \cite{HaT}).
In particular, it is known from early work of Erd\H os \cite{E},
that for each $\eps>0$ there is some $\eta>0$, such that
if $S$ is the set of integers in $(N,N^{1+\eta}]$, then
$\delta^+(S)\ge1-\eps$ for all large $N$.  In fact, we
almost have an asymptotic estimate for $1-\delta^+(S)$ for
such a set $S$:
Among other results, it is shown in
Theorem 1 of \cite{F} that for $0<\eta<1/2$ and
$N\ge 2^{1/\eta}$, $\delta^+(S)$ is between 
$1-c_1\eta^{\theta}(\log1/\eta)^{-3/2}$ and $1-c_2\eta^{\theta}(\log1/\eta)^{-3/2}$, 
where $c_1,c_2$ are positive absolute constants and
where $\theta=1-(1+\log\log 2)/\log 2 =0.08607\ldots$.

In the above example with $\eta>0$ fixed, we have
$\delta^-(S)\le\alpha(S)=(1+o(1))N^{-\eta}=o(1)$,
while for large $N$, $\delta^+(S)$ is bounded away from 0.
If the residue classes are chosen randomly, should we expect the density 
of the residual set to be closer to $\delta^-(S)$, $\alpha(S)$, or $\delta^+(S)$?  
We show in Sections~\ref{near} and \ref{normal} 
that for any finite integer set $S$, the average (and typical) 
case has residual density close to~$\alpha(S)$.

Finally we mention a problem we have not been able to settle.
Is it true that for each positive number $B$, there are positive
numbers $\Delta_B$, $N_B$, such that if $S$ is a finite set of positive
integers greater than $N_B$ with reciprocal sum at most $B$, 
then $\delta^-(S)\ge \Delta_B$?  If this holds it would imply
each of Conjectures 1, 2, and 3.
For more problems and results concerning covering systems, the
reader is directed to \cite{PS} and \cite{S}.

\medskip
\noindent{\bf Acknowledgements.}  We would like to thank the
referee for insightful comments regarding our inequality
\eqref{keylemma}, and in particular for the content of Remark~\ref{lllrem}
in the next section.  We also would like to
thank G. Tenenbaum for informing us of the theorem of Rogers
mentioned above.

\section{A basic lemma and Haight's theorem}

To set some notation, we shall always have $n$ a positive integer,
with $P(n)=P^+(n)$ the largest prime factor of $n$ for $n>1$ and $P(1)=0$.  
We shall also let $P^-(n)$ denote the least prime factor of $n$ when
$n>1$, and $P^-(1)=+\infty$.
The letter $p$ will always represent a prime variable. 
We use $N,K,Q$ to represent real numbers, usually large.
We use the Vinogradov notation $\ll$ from analytic number theory,
so that $A\ll B$ is the same as $A=O(B)$, but is cleaner to use
in a chain of inequalities.  In addition, $A\gg B$ is the
same as $B\ll A$.  All constants implied by this
notation are absolute and bounds for them are computable in principle.
If $S$ is a multiset, and we have some product or sum with $n\in S$,
it is expected that $n$ is repeated as many times in the product or
sum as it appears in $S$.

Let $C$ be a finite set of ordered pairs of positive integers $(n,r)$,
which we interpret as a set of residue classes $r \pmod n$.
We say such a set is a {\it residue system}.
Let $S=S(C)$ be the multiset of the moduli $n$ appearing in $C$.
The number of times an integer $n$ appears in $S$ we call the
\emph{multiplicity} of $n$.
By $R(C)$ we denote the set of integers {\it not} congruent to
$r\pmod n$ for any $(n,r)\in C$.  Since $R(C)$ is a union of residue classes
modulo the least common multiple of the members of $S(C)$, it follows
that $R(C)$ possesses a (rational) asymptotic density, which we denote by
$\delta(C)$.
If $C = \{ (n_1,r_1), \ldots, (n_l,r_l)\}$, then we set
\[
\alpha(C)=\prod_{n\in S(C)}\(1-\frac1{n}\) = \prod_{j=1}^l\(1-\frac1{n_j}\),
\quad
\beta(C)=\sum_{\substack{i<j \\ \gcd(n_i,n_j)>1}} \frac1{n_i n_j}.
\]
Note that $\alpha(C)$ depends only on $S(C)$, so it is notationally
consistent with $\alpha(S)$ from Section~\ref{intro}.


\begin{lem}\label{basiclem}
For any residue system $C$, we have $\delta(C)\ge\alpha(C)-\beta(C)$.
\end{lem}

\begin{proof}  Let $\alpha=\alpha(C)$ and $\beta=\beta(C)$.
We use induction on $l$. If $l=1$, then $\beta=0$ and the
statement is trivial.  Let $l>1$; we will describe an induction
step from $l-1$ to $l$.  We denote $C'=\{(n_1,r_1),\ldots,(n_{l-1},r_{l-1})\}$,
~$C''=\{(n_j,r_j):\,j<l,~\gcd(n_j,n_l)=1\}$,
\[
\a'=\a(C') = \prod_{j=1}^{l-1}\left(1-\frac1{n_j}\right)~\hbox{ and }~
\b'=\b(C') = \sum_{\substack{ i<j\le l-1,\\ \gcd(n_i,n_j)>1}} \frac1{n_in_j}.
\]

By the induction supposition,
\be\label{lem1a}
\delta(C')\ge\alpha'-\beta'.
\ee
Also,
\be\label{lem1b}
\delta(C'')\le\delta(C')+\sum_{n_j\in S(C')\setminus S(C'')}\frac1{n_j}
=\delta(C')+\sum_{\substack{j<l \\ \gcd(n_j,n_l)>1}}\frac1 {n_j}.
\ee
The density of integers covered by the residue class $r_l \pmod{n_l}$
but not covered by $r_j \pmod{n_j}$ for every $n_j\in S(C'')$ is equal to
$\delta(C'')/n_l$.  Therefore,
\begin{align*}
\delta(C')-\delta(C)
&=\hbox{density}\{n\equiv r_l~({\rm mod}~n_l):n\in R(C')\}\\
&\le\hbox{density}\{n\equiv r_l~({\rm mod}~n_l):n\in R(C'')\}
=\delta(C'')/n_l,
\end{align*}
so that, by \eqref{lem1a} and \eqref{lem1b},
\begin{align*}
\delta(C)&\ge\delta(C')-\Bigg(\delta(C')
      +\sum_{\substack{j<l \\
\gcd(n_j,n_l)>1}}\frac1{n_j}\Bigg)\frac1{n_l}\\[5pt]
&=\left(1-\frac1{n_l}\right)\delta(C')
      -\sum_{\substack{j<l \\ \gcd(n_j,n_l)>1}}\frac1{n_jn_l}\\[5pt]
&\ge\left(1-\frac1{n_l}\right)(\alpha'-\beta')-(\beta-\beta')
      \ge\left(1-\frac1{n_l}\right)\alpha'-\beta=\alpha-\beta.
\end{align*}
This completes the proof of the lemma.
\end{proof}

\begin{remark}
The proof of Lemma \ref{basiclem} actually gives the better bound
\[
\delta(C) \ge \a(C) - \sum_{\substack{ i<j \\ \gcd(n_i,n_j)>1 }}
\frac{1}{n_i n_j} \prod_{u>j} \bigg(1-\dfrac{1}{n_{u}}\bigg).
\]
\end{remark}
\begin{remark}
\label{lllrem}
The referee has pointed out to us that Lemma \ref{basiclem} can
be formulated in a more general way involving a finite number
of events in a probability space.  In particular suppose that
$E_1,E_2,\dots,E_l$ are events in a probability space with the property 
that if $E_i$ is independent individually of the events 
$E_{j_1},E_{j_2},\dots E_{j_k}$, then it is independent of every
event in the sigma algebra generated by $E_{j_1},E_{j_2},\dots,E_{j_k}$.
Then
\begin{equation} \label{lovasz}
{\rm P}\left({\textstyle\bigcap_{i=1}^l\overline E_i}\right)\ge
\prod_{i=1}^l{\rm P}(\overline{E_i})
-\sum_{\substack{1\le i<j\le l\\ E_i,E_j~{\rm dependent}}}{\rm P}(E_i){\rm P}(E_j).
\end{equation}
We can retrieve Lemma \ref{basiclem} from this statement if we let
$E_i$ be the event that an integer $n$ is in the residue class $r_i\pmod{n_i}$.
Indeed, $E_i$ is independent of $E_j$ if and only if $n_i$ and $n_j$
are coprime.  The extra condition involving the sigma algebra is easily
seen to hold (and was used strongly in our proof).  The proof of 
\eqref{lovasz} is the same as that of Lemma \ref{basiclem}, 
namely an induction on $l$.  This result bears a resemblance to the 
Lov\'asz Local Lemma
(for example, see \cite{AS}), and may be stronger than it in some situations.  
\end{remark}


There is a very interesting negative result of Haight \cite{H}.
As in the introduction, we say an integer $H$ is covering
if there is a covering system with the moduli being the (distinct)
divisors of $H$ that are larger than 1.
It is shown in \cite{H} that there exist integers $H$
that are {\em not} covering, yet $\sum_{d\mid H}1/d=\sigma(H)/H$ is arbitrarily
large.
Although Haight's theorem follows directly from Theorem A
(by taking $K$ fixed, $N$ large
and $H$ to be the product of the primes in $(N,N^K]$),
Lemma \ref{basiclem} by itself leads to a
new (and short) proof of a stronger version of Haight's result:
\medskip

\begin{thm}\label{Haight}
There is an infinite set of positive integers $H$ with
\[
\sigma(H)/H=(\log\log H)^{1/2}+O(\log\log\log H),
\]
such that for
any residue system $C$ with $S(C)=\{d:d>1,~d\mid H\}$, we have
\[
\delta(C)\ge(1+o(1))\alpha(C).
\]
In particular, for large $H$ in this set, no such $C$
can have $\delta(C)=0$.
\end{thm}

\begin{proof}
Let $N$ be a large parameter, and let
\[
H=\prod_{e^{\sqrt{\log N}}\log N<p\le N}p.
\]
Then
\[
\log\sum_{d\mid H}\frac1d=\sum_{e^{\sqrt{\log N}}\log N<p\le N} \(\frac1p
+O\pfrac{1}{p^2}\)
=\frac12 \log\log N -\frac{\log\log N+O(1)}{\sqrt{\log N}},
\]
by Mertens' theorem.  Thus, as $\log H =(1+o(1))N$ by the prime number
theorem, we have
\[
\frac{\sigma(H)}{H} =
\sum_{d\mid H}\frac1d=(\log N)^{1/2}-\log\log N+O(1)=(\log\log H)^{1/2}
+O(\log\log\log H).
\]
Let $C$ be a residue system with $S(C)=\{d:d>1,~d\mid H\}$.  We have
\[
\log\alpha(C)=\sum_{d\in S(C)}\log(1-1/d)=-\sum_{d\in 
S(C)}1/d+O\(\exp(-\sqrt{\log N})\),
\]
so that
\[
\alpha(C)=\exp\(-\sqrt{\log N}+O(1)\)\log N.
\]
Also,
\[
\beta(C)\le\sum_{d>1}\sum_{\substack{d_1,d_2\in S(C)\\ d\mid 
d_1,~d\mid d_2}}\frac1{d_1d_2}
\le\sum_{d\mid H,~d>1}\frac1{d^2}\sum_{\substack{d_1\mid H\\ d_2\mid 
H}}\frac1{d_1d_2}
\ll\log N\sum_{d\mid H,~d>1}\frac{1}{d^2}.
\]
Further,
\[
\sum_{d\mid H,~d>1}\frac1{d^2}\le
\sum_{d>e^{\sqrt{\log N}}\log N}\frac1{d^2}
\ll\exp\(-\sqrt{\log N}\)(\log N)^{-1}.
\]
Thus,
\[
\beta(C)\ll\exp\(-\sqrt{\log N}\)=o(\alpha(C))
\]
and the theorem follows from Lemma \ref{basiclem}.
\end{proof}

\begin{remark}
An examination of our proof shows that we have a more general result.
Let $\cal H$ be the set of integers $H$ which have no prime factors
below $\exp(\sqrt{\log\log H})\log\log H$.  As $H\to\infty$ in $\cal H$
we have for any residue system $C$ with $S(C)=\{d:d>1,~d\mid H\}$ that
$\delta(C)\ge(1+o(1))\alpha(C)$.  In particular, at most finitely many
integers $H\in\cal H$ are covering.

We also remark that the proof gives the following result.  Say that
a positive integer $H$ is $s$-covering, if for each $d\mid H$ with $d>1$
there are $s$ integers $r_{d,1},\dots,r_{d,s}$ such that the union of
the residue classes $r_{d,i}~({\rm mod}~d)$ for $i=1,\dots,s$, and
$d\mid H$ with $d>1$ is $\Z$.  Then for each fixed $\eps>0$ there
are values of $H$ where $\sigma(H)/H$ is arbitrarily large,
yet $H$ is not $s$-covering with $s=[(\log\log H)^{1-\eps}]$.
Indeed, take $H$ to be the product of the primes in
$\left(\exp\((\log N)^{1-\eps/3}\),N\right]$ and follow the same proof.
This too strengthens a result in \cite{H}.
\end{remark}

\section{The smooth number decomposition}

The relative ease of using Lemma \ref{basiclem} in the proof of Haight's
theorem is due to the fact that the moduli that we produce for the
proof have no small prime factors, so that it is easy to bound
the sum for $\beta(C)$.  In going over to more general cases
it is clear we have to introduce other tools.  For example, if
$S(C)$ is the set of all integers in the interval $(N,KN]$, then
the sum for $\beta(C)$ tends to infinity with $K$, while the
expression for $\alpha$ is always less than 1.  Thus, the lemma
would say that the residual set of integers not covered has
density bounded below by a negative quantity tending to $-\infty$.
This is clearly not useful!  To rectify this situation, we
choose a paramter $Q$ and factor each modulus $n$ as 
$n_{\underline Q}n_{\overline{Q}}$, where $n_{\underline Q}$ is 
the largest divisor of $n$ composed solely of
primes in $[1,Q]$, and $n_{\overline{Q}}=n/n_{\underline Q}$.  We then
find a way to decompose our system $C$ based on these factorizations, 
and use Lemma~\ref{basiclem} on the parts corresponding to the
numbers $n_{\overline{Q}}$ which have no small prime factors.  

To set some terminology, for a number $Q\ge1$,
we say a positive integer $n$ is $Q$-smooth if $P(n)\le Q$.
Thus, $n_{\underline Q}$ is the largest $Q$-smooth divisor of $n$.


\begin{lem}\label{decomposition}
Let $C$ be an arbitrary residue system.  
Let $Q\ge2$ be arbitrary, and set
\[
M=\lcm\{n_{\underline Q}:n\in S(C)\}.
\]
For $0\le h\le M-1$, let $C_h$ be the set
\[
C_h = \left\{ \(n_{\overline{Q}},r  \) : (n,r) \in C,~r  \equiv h
{\hskip -8pt} \pmod{n_{\underline Q}} \right\}.
\]
Then
\[
\delta(C)=\frac1M\sum_{h=0}^{M-1} \delta(C_h).
\]
\end{lem}

\begin{proof}
Fix $h$ so that $0\le h\le M-1$.
For $(n,r) \in C$, the simultaneous congruences
\[
x\equiv r\hbox{ (mod }n\hbox{)}, \qquad x\equiv h\hbox{ (mod }M\hbox{)}
\]
have a solution if and only if $r \equiv h \pmod{n_{\underline Q}}$, 
since $n_{\underline Q}=\gcd(n,M)$,
in which case the system is equivalent to the system
\[
x\equiv r\kern-8pt\pmod{n_{\overline{Q}}},\qquad x\equiv h\kern-8pt\pmod{M}.
\]
Thus, $R(C_h)\cap(h~{\rm mod}~M)=R(C)\cap(h~{\rm mod}~M)$.
Observe that all elements 
$n_{\overline{Q}}$  
of $S(C_h)$ are coprime to $M$.
Thus, the proportion of the numbers in $R(C_h)$ in the class
$h$ modulo $M$ is equal to $\delta(C_h)$.
Hence, the density of $R(C)\cap(h~{\rm mod}~M)$ is $\delta(C_h)/M$
and the result follows.
\end{proof}


We now take advantage of the fact that the prime factors of
a number $n_{\overline{Q}}$ are all larger than $Q$ to allow
us to get a reasonable upper bound for the quantities $\beta(C_h)$.
The proof is similar to that in Theorem~\ref{Haight}.

\begin{lem}\label{sumbeta}
Let $K>1$, and suppose 
$C$ is a residue system with $S(C)$ consisting
of integers in the interval
$(N,KN]$, each with multiplicity at most $s$.
Suppose $Q\ge 2$, and define 
$M$ and $C_h$ as in Lemma~\ref{decomposition}.  Then
\be\label{sumb}
\frac1M \sum_{h=0}^{M-1} \beta(C_h) \ll  \frac{s^2\log^{2} (QK)}{Q}.
\ee
\end{lem}
\begin{proof}
For $m|M$, let $S_m$ be the set of distinct numbers 
$n_{\overline{Q}}=n/\gcd(n,M)$, 
where $n\in S(C)$ and 
$n_{\underline Q}=\gcd(n,M)=m$.
For $m,m'\mid M$, let
$$
F(r,m,r',m') = \# \{ 0\le h\le M-1 : h\equiv r \hbox{ (mod } 
m\hbox{)},~h\equiv r' \hbox{ (mod }m'\hbox{)} \}.
$$
Then
$$
\frac1M \sum_{h=0}^{M-1} \beta(C_h) \le \frac1M \sum_{\substack{m|M\\m'|M}}
      \sum_{\substack{ n\in S_{m} \\ n'\in S_{m'} \\  \gcd(n,n')>1}}
      \frac{1}{nn'} \sum_{\substack{(nm,r)\in C \\
      (n'm',r')\in C}} F(r,m,r',m').
$$
Since $F(r,m,r',m')$ is either $0$ or $M/\lcm[m,m']$, the 
inner sum is at most
$$
s^2 \frac{M}{\lcm[m,m']}.
$$
Next,
\benn\begin{split}
\sum_{\substack{ n\in S_{m}\\n'\in S_{m'}\\\gcd(n,n')>1}}
\frac1{nn'}
&\le\sum_{p>Q}\sum_{\substack{ n\in S_{m}\\n'\in S_{m'}\\p|n,~p|n'}}
\frac1{nn'}\\[5 pt]
&= \sum_{p>Q}\Biggl(\sum_{\substack{ N/m<n\le KN/m\\p|n,~P^-(n)>Q}}
\frac{1}{n}\Biggr)\Biggl(\sum_{\substack{ N/m'<n'\le KN/m'\\
p|n',~P^-(n')>Q}} \frac{1}{n'}\Biggr).
\end{split}
\eenn
By standard sieve methods (e.g., Theorem 3.3 of \cite{HR}),
uniformly in $x\ge 2$, $z\ge 2$, the number of integers $\le x$ which
have no prime factor $\le z$ is $\ll x/\log z + 1$.
By partial summation,
\[
\sum_{\substack{ N/m<n\le KN/m\\p|n,~P^-(n)>Q}}\frac1{n} =
\frac{1}{p} \sum_{\substack{\frac{N}{pm} < t \le \frac{KN}{pm} \\
        P^-(t)>Q} } \frac{1}{t}
\ll \frac1p\( \frac{\log K}{\log Q}+1 \)
=\frac{\log(QK)}{p\log Q}
\]
and similarly with $m',n'$ replacing $m,n$.  
We have the estimate $\sum_{p>Q}p^{-2}\ll 1/(Q\log Q)$, 
which follows
from the prime number theorem and partial summation.  Thus,
\[
\sum_{\substack{ n\in S_{m}\\n'\in S_{m'}\\\gcd(n,n')>1}}
\frac1{nn'}
\ll\frac{\log^2(QK)}{Q\log^3Q},
\]
so that
\[
\frac1M \sum_{h=0}^{M-1} \beta(C_h) \ll
\frac{s^2 \log^2 (QK)}{Q\log^3 Q} \sum_{\substack{m|M\\m'|M}}
      \frac1{\lcm[m,m']}
= \frac{s^2 \log^2 (QK) }{Q\log^3 Q} \sum_{u|M}
      \sum_{\substack{m|M,~m' | M \\ \lcm[m,m']=u}} u^{-1}.
\]
With $\tau(n)$ denoting the number of natural divisors of $n$, the double sum
is equal to
\[
\sum_{u|M} u^{-1} \tau(u^2)
\le\prod_{p|M}\(1+\frac3p +\frac{5}{p^2}+\cdots\)=\prod_{p|M}\frac{1+1/p}{(1-1/p)^2}
\le\prod_{p\le Q}\frac{1+1/p}{(1-1/p)^2}
\ll \log^3 Q,
\]
and this completes the proof.
\end{proof}

%
%

To complement Lemma~\ref{sumbeta}, we would like
a lower bound for the sum of the $\alpha(C_h)$.
Key to this estimate will be those moduli in $S(C)$
which are $Q$-smooth.  If the residue classes corresponding to these
moduli do not cover everything, we are able to get
a respectable lower bound for the sum of the $\alpha(C_h)$.

\begin{lem}\label{sumalpha}
Suppose that $C$ is a
residue system, $Q\ge 2$, 
and define $M$ and $C_h$ as in Lemma \ref{decomposition}.  Also let 
$C'= \{ (n,r) \in C: n | M \}=\{(n,r)\in C:P(n)\le Q\}$
and suppose  $\delta(C') > 0$. Then
\[
\frac1M\sum_{h=0}^{M-1} \alpha(C_h)
\ge \left(\alpha(C)\right)^{(1+1/Q)/\delta(C')}.
\]
\end{lem}

\begin{proof}
Note that $1\in S(C_h)$ if and only if there is a pair $(n,r)\in C'$ 
with $h \equiv r \pmod{n}$.  Let
\[
\mathcal M' = \{ 0 \le h\le M-1 : 1\not\in
S(C_h)\}, \quad M' = |\mathcal M'|.
\]
Then
\be\label{ratio}
\frac{M'}{M} = \delta(C').
\ee
The hypothesis $\delta(C')>0$ thus implies that $M'>0$.
Observe that $1\in S(C_h)$ implies $\alpha(C_h) =0$.
By the inequality of the arithmetic and geometric means,
\benn\begin{split}
\frac1M \sum_{h=0}^{M-1} \alpha(C_h) &=
      \frac1M \sum_{h\in\mathcal M'} \alpha (C_h) \ge \frac {M'}M
      \biggl(\prod_{h\in \mathcal M'} \alpha(C_h) \biggr)^{1/M'} \\[5pt]
&=\frac {M'}M \bigg( \prod_{h\in \mathcal M'} \prod_{n'\in S(C_h)}
       \(1- \frac{1}{n'}\)  \bigg)^{1/M'}.
\end{split}
\eenn
Since $\log(1-1/k)>-\frac{1}{k}(1+\frac1k)$ for $k\ge2$ and since
each $n' > Q$, we have
\[
1 - \frac{1}{n'} > \exp \big({\hskip -2pt}- \lam / n' \big), \quad
\text{ where } \lam=1+1/Q.
\]
Thus,
\benn
\begin{split}
\frac1M \sum_{h=0}^{M-1} \alpha(C_h) &\ge \frac{M'}{M} \exp\(
-\frac{\lam}{M'} \sum_{h\in\mathcal M'} \sum_{n'\in S(C_h)}
\frac{1}{n'} \) \\[5pt]
&\ge \frac{M'}{M} \exp \( \frac{\lam(M-M')}{M'} -
\frac{\lam}{M'} \sum_{h=0}^{M-1} \sum_{n'\in S(C_h)}
\frac{1}{n'} \),
\end{split}
\eenn
where the last inequality uses that $1\in S(C_h)$ for $h\not\in \mathcal M'$.

Each pair $(n,r)\in C$
maps to those $C_h$ with $h\equiv r\kern-5pt\pmod{n_{\underline Q}}$, so it produces
the pair $(n_{\overline{Q}},r)$ in exactly $M/n_{\underline Q}$ sets $C_h$ for
$h\in[0,M-1]$.  
We thus have
\[
\sum_{h=0}^{M-1} \sum_{n'\in S(C_h)} \frac{1}{n'}
\le \sum_{n\in S(C)}\frac{M}{n_{\underline Q}}\cdot\frac1{n_{\overline Q}}
= M \sum_{n\in S(C)} \frac{1}{n}.
\]
(Note that the inequality holds since several pairs in $C$ may map to the same
pair in some $C_h$, where they would be counted just once.)  Thus,
\[
\frac1M\sum_{h=0}^{M-1}\alpha(C_h)\ge\frac{M'}{M}\exp
\left(\frac{\lambda(M-M')}{M'}-\frac{\lambda M}{M'}\sum_{n\in 
S(C)}\frac1n\right).
\]
Also, $(M'/M) \exp\big( (M-M')/M' \big) \ge 1$.  Thus,
\[
\frac1M \sum_{h=0}^{M-1} \alpha(C_h) \ge \exp\( -\frac{\lam M}{M'}
\sum_{n\in S(C)} \frac{1}{n}\) \ge \bigl( \a(C) \bigr) ^{\lam M/M'}.
\]
The lemma follows by \eqref{ratio}.
\end{proof}

We now combine our lemmas into one easily-applied statement.

\begin{lem}\label{combinedlemma}
Suppose $K>1$, $N$ is a positive integer, and $C$ is a residue system with
$S(C)$ consisting of integers in $(N,KN]$, each with multiplicity at most $s$.
Let $Q\ge2$, and as in Lemma~\ref{sumalpha}, let $C'=\{(n,r)\in C:P(n)\le Q\}$.
If $\delta(C')>0$, then
\[
\delta(C)\ge\alpha(C)^{(1+1/Q)/\delta(C')}+O\(\frac{s^2\log^{2} (QK)}{Q}\),
\]
where the implied constant is uniform in all parameters.
\end{lem}

\begin{proof}
Define $M$ and $C_h$ as in Lemma \ref{decomposition}.
By Lemmas~\ref{basiclem}, \ref{decomposition}, \ref{sumbeta}, and 
\ref{sumalpha}, we have
\benn
\begin{split}
\delta(C)&=\frac1M \sum_{h=0}^{M-1} \delta(C_h)
\ge\frac1M\sum_{h=0}^{M-1} \alpha(C_h)
-\frac1M\sum_{h=0}^{M-1}\beta(C_h)\\
&\ge \a(C)^{(1+1/Q)/\delta(C')} +
O\(\frac{s^2 \log^{2} (QK)}{Q} \).
\end{split}
\eenn
Thus, we have the lemma.
\end{proof}

%
%
\section{Lower bounds on $\delta(C)$}\label{lowerbounds}
%
%

In this section we prove stronger versions of Theorems A and B.
We begin with a useful lemma about smooth numbers.

\begin{lem}\label{smoothsum}
Suppose $Q\ge 2$ and $Q < N \le \exp (\sqrt{Q})$.  Then
$$
\sum_{\substack{n>N \\ P(n) \le Q}} \frac{1}{n} \ll (\log Q) e^{-u\log
   u}, \quad \text{where } u=\frac{\log N}{\log Q}.
$$
\end{lem}

\begin{proof}
We use standard upper-bound estimates for
the distribution of smooth numbers:
The number of $Q$-smooth numbers
at most $t$ is $\ll t/u_t^{u_t}$, where $u_t=\log t/\log Q$, provided
$Q\le t\le \exp\(Q^{1-\eps}\)$ (\cite{HT}, Theorem 1.2 and Corollary
2.3).   Further, for $t>\exp\(6\sqrt{Q}\)$,
the $Q$-smooth numbers are distributed more sparsely than the squares.
We thus have
\benn
\begin{split}
\sum_{\substack{n>N \\ P(n)\le Q}} \frac{1}{n}
&=\int_N^\infty \frac1{t^2}\sum_{\substack{N<n\le t\\ P(n)\le Q}}1\,dt
\le\sum_{0\le i\le 
10\sqrt{Q}}\int_{NQ^i}^{NQ^{i+1}}\frac1{t^2}\sum_{\substack{n\le t\\ 
P(n)\le Q}}1\,dt
+\int_{\exp(6\sqrt{Q})}^\infty \frac1{t^2}\sum_{\substack{n\le t\\ 
P(n)\le Q}}1\,dt\\
&\ll \sum_{i\ge0}\frac{\log Q}{(u+i)^{u+i}}
+\int_{\exp(6\sqrt{Q})}^\infty \frac1{t^{3/2}}\,dt
\ll\frac{\log Q}{u^u},
\end{split}
\eenn
implying the lemma.
\end{proof}

Let
$$
L(N,s) =\exp \bigg( \log N \frac{\log\log(s\log N)}{\log(s\log N)} \bigg).
$$


\begin{thm}\label{thmsmallsum}
Suppose $0<b<\frac12$, $0<c<\frac13(1-4b^2)$ and
let $N$ be sufficiently large, depending on the choice of $b$ and $c$.
Suppose $C$ is a
residue system with $S(C)$ consisting of integers $n>N$,
each having multiplicity at most $s$, where $s \le \exp\big(b\sqrt{\log N \log
   \log N}  \big)$, and such that
\be\label{largesum}
\sum_{n\in S(C)}\frac1n\le c \log L(N,s).
\ee
Then $\delta(C)>0$.
\end{thm}

\begin{proof}
Throughout we assume that $N$ is sufficiently large, depending only on $b$ and $c$.
Let $\lambda = \frac13(1-4b^2)$ and put $\eps=\frac{1}{20}(\lambda-c)$.
First, we have
$$
-\log \a(C) \le \sum_{n\in S(C)} \(\frac{1}{n}+\frac{1}{n^2}\)
\le \bigg(  1 + \dfrac{1}{N}  \bigg) \sum_{n\in S(C)} \frac{1}{n}
\le (c+\eps)\log L(N,s) = G,
$$
say.  Define
\be\label{Qdef}
Q_0=L(N,s)^{1-\eps}, \quad Q_{j} = \exp \( Q_{j-1}^{\lambda+\eps} \)
\quad (j\ge 1)
\ee
and
$$
K_j = \exp \( Q_{j-1}^{\lambda+2\eps} \) \qquad (j\ge 1).
$$
Let
$$
C_j = \{ (n,r)\in C : P(n) \le Q_j \}.
$$
Also, define
\be\label{deltadef}
\del_0 = 1 - \eps, \quad \del_j = e^{-1-G(1+1/Q_0)/\del_{j-1}} \qquad
(j\ge 1),
\ee
where $G$ is defined above.  Since $C$ is finite and $Q_j$ tends to 
infinity with $j$,
it follows that $C=C_j$ for large $j$.  Thus, the theorem will follow 
if we show that
\be\label{deltaj}
\delta(C_j) \ge \del_j \qquad (j\ge 0).
\ee

First, by Lemma \ref{smoothsum},
$$
1-\delta(C_0) \le s
\sum_{\substack{n>N \\ P(n) \le Q_0}} \frac{1}{n}
\ll s (\log Q_{0} ) e^{-u\log u}, \quad \text{where } u = \frac{\log 
N}{\log Q_0}.
$$
By the definition of $Q_0$, we have
\[
u = \dfrac{\log (s \log N)}{(1-\eps) \log \log (s \log N)}
\]
so that
$\log u\ge(1-\eps)\log\log(s\log N)$ and
$u\log u \ge \log(s\log N)$.
Hence, $\delta(C_0) \ge \del_0$.

Next, suppose $j\ge 1$ and $\del(C_{j-1}) \ge \del_{j-1}$.
Let $s_0 =\exp\big( b\sqrt{\log N \log \log N} \big)$ and observe that
for $N$ large and $s\le s_0$, we have
\[
\dfrac{\log\log(s \log N)}{\log(s \log N)}
\ge \dfrac{\log\log s_0}{\log (s_0 \log N)}
\ge \dfrac{\log \log N}{2 \log (s_0 \log N)}
\ge \dfrac{(1-\eps) \log \log N}{2b \sqrt{ \log N \log\log N}}.
\]
Therefore,
\[
s^{2} \le \exp\big( 2b\sqrt{\log N \log \log N} \big)
\le L(N,s)^{4b^2/(1-\eps)}
= Q_{0}^{4b^2/(1-\eps)^2} \le Q_{0}^{4b^2(1+3\eps)}\le Q_{j-1}^{4b^2(1+3\eps)}.
\]
Let
$$
C_j' = \{ (n,r)\in C_j : n\le K_j \}, \quad C_j''  = \{ (n,r)\in C_j : n >
K_j \}.
$$
Observe that
\[
\delta \big(  \{ (n,r)\in C'_j : P(n) \le Q_{j-1} \}  \big)
\ge \delta(C_{j-1}) \ge \delta_{j-1}
\quad \text{ and } \quad
\alpha(C'_{j}) \ge \alpha(C) \ge e^{-G}.
\]
By Lemma \ref{combinedlemma}
with $Q=Q_{j-1}$ and $K=K_j/N$, there is an absolute constant $D$ such that
\be\label{induct1}
\begin{split}
\del(C_j') &\ge \a(C_j')^{(1+1/Q_{j-1})/\del_{j-1}} - D
   \frac{s^{2}\log^2 (Q_{j-1}K_j/N)}{Q_{j-1}} \\
&\ge e^{-G(1+1/Q_0)/\del_{j-1}} - Q_{j-1}^{-1+4b^2(1+3\eps)+2\lambda+5\eps}
\ge 2\delta_j- Q_{j-1}^{-\lambda + 8\eps}.
\end{split}
\ee
Also, by Lemma \ref{smoothsum},
\be\label{induct2}
1-\del(C_j'') \le s
\sum_{\substack{n>K_j \\ P(n) \le Q_{j}}} \frac{1}{n} \ll
s (\log Q_{j}) e^{-u_j\log u_j},
\ee
where
$$
u_j=\frac{\log K_j}{\log Q_j} = Q_{j-1}^{\eps}.
$$
Thus, $1-\del(C_j'')\le Q_{j-1}^{-1}$.
Together with \eqref{induct1}, this implies
\be\label{induct3}
\del(C_j) \ge \del(C_j') - (1-\del(C_j'')) \ge 2\del_j -
Q_{j-1}^{-\lambda+9\eps}.
\ee

To complete the proof of \eqref{deltaj} and the theorem, it suffices
  to prove that
\be\label{induct4}
Q_{j-1}^{-\lambda+9\eps} \le \del_j \qquad (j\ge 1).
\ee
First,
$$
Q_0^{-\lambda+9\eps} = L(N,s)^{(-\lambda+9\eps)(1-\eps)} \le
L(N,s)^{-\lambda+10\eps}=L(N,s)^{-c-10\eps},
$$
while
$$
\del_1 \ge e^{-1-G(1+1/Q_0)(1+1.1\eps)} \ge L(N,s)^{-c-2\eps}.
$$
This proves \eqref{induct4} when $j=1$.
Suppose \eqref{induct4} holds for some $j\ge1$.  Since $G\le \log
Q_0$, we have
\[
-\log\del_{j+1}=1+G(1+1/Q_0)/\del_j \le 2G Q_{j-1}^{\lambda-9\eps} <
  Q_{j-1}^{\lambda-8\eps}.
\]
And, by \eqref{Qdef}, we have $-\log(Q_j^{-\lambda+9\eps})\gg Q_{j-1}^{\lambda+\eps}$.
Thus, $Q_{j}^{-\lambda+9\eps} \le \del_{j+1}$ and by induction
\eqref{induct4} holds for all $j$.
This completes our proof.
\end{proof}

Theorem \ref{thmsmallsum} implies Theorem A of the introduction by 
setting $s = 1$.
Observe that the bound on the sum in Theorem \ref{thmsmallsum} given 
in (\ref{largesum})
decreases as $s$ increases.
If one is interested in a result similar to Theorem A but with an 
emphasis on allowing the multiplicity
of the moduli to be large, one may take $b$ arbitrarily close to $1/2$
in Theorem \ref{thmsmallsum}.

Theorem \ref{thmsmallsum} should be compared with Theorem 
\ref{coveringexample} of the
next section which shows that coverings, even exact coverings with 
squarefree moduli, exist when we
allow the multiplicity of the moduli to be of size
$\exp\big( \sqrt{\log N \log \log N} \big)$.

We can also consider the case that $S(C)$ consists of integers
from $(N,KN]$ with multiplicities at most $s \le
\exp\big(b \sqrt{\log N \log \log N} \big)$,
where $b < \sqrt{3\eps/4}$.
If $0<\eps<1/3$, $N$ is large, and $K = L(N,s)^{(1/3-\eps)/s}$, then
Theorem \ref{thmsmallsum} implies that $\delta(C)>0$.
By a different argument, we can extend the range of $K$ a bit.

\begin{thm}\label{thmpos}
Suppose $0<\eps<(1-\log 2)^{-1}$, $b<\frac12\sqrt{(1-\log 2)\eps}$
and $N$ is sufficiently large, depending on the choice of $\eps$ and $b$.
Suppose that $C$ is a residue system with $S(C)$ consisting of
integers from $(N,KN]$
with multiplicity at most $s$, where
$s \le \exp\big(b \sqrt{ \log N \log \log N}  \big)$
and $K =  L(N,s)^{( (1-\log 2)^{-1}-\eps)/s}$.
Then $\delta(C)>0$.
\end{thm}

Note that for $s\ge \log N$, $K=1+o(1)$.
Before proving Theorem \ref{thmpos}, we present a lemma.

\begin{lem}\label{lem2pt1}
Suppose $s$ is a positive integer and $C$ is a residue system with 
$S(C)$ consisting
of integers from $(1,B]$ with multiplicity at most $s$.  Let
\[
C_0=\{(n,r)\in C:P(n)\le\sqrt{s B}\}.
\]
If $\delta(C_0)>0$, then $\delta(C)>0$.
\end{lem}

\begin{proof}
Suppose that $\delta(C_0) > 0$.
Denote by $P$ the product of all primes in $(\sqrt{s B},B]$, and
let $L$ be the least common multiple of the elements of $S(C_0)$.

Let $p$ be a prime divisor of $P$.
Since $p > \sqrt{s B}$ implies $s B/p < p$, there are at most $p-1$ 
multiples of
$p$ in the multiset $S(C)$.  Call them $m_{1}, \dots, m_{t}$, and let
$r_1,\ldots, r_t$ be the corresponding residue classes.   Then there is a
choice for $b = b(p) \in \{ 0, 1, \dots, p-1 \}$ such that each integer
satisfying $x \equiv b \pmod{p}$ is
not covered by (i.e., does not satisfy) any of the congruences
$x\equiv r_j \pmod{m_{j}}$ with $1 \le j \le t$.

By assumption, there is a residue class $a~{\rm mod}~L$ contained in $R(C_0)$.
Let $A$ be a solution to the Chinese remainder system $A\equiv a\pmod L$ and
$A\equiv b(p)\pmod p$ for each prime $p$ dividing $P$.  Then not only 
do we have
$A\not\equiv r\pmod n$ for each $(n,r)\in C_0$, we also have for each 
prime $p\mid P$
and $(n,r)\in C$ with $p\mid n$, that $A\not\equiv r\pmod n$.  Since 
this exhausts
the pairs $(n,r)\in C$, we have $A\in R(C)$, so we have the lemma.
\end{proof}

\begin{proof}[Proof of Theorem \ref{thmpos}]
We may suppose that $\eps>0$ is sufficiently small and $K\ge 2$.
Let $C_0$ be as in Lemma \ref{lem2pt1} where we take $B=KN$.  Then
\begin{align*}
\sum_{n\in S(C_0)}\frac1n
&\le s \sum_{\substack{N<n\le KN\\ P(n)\le\sqrt{s KN}}}\frac1n
= s \sum_{N<n\le KN}\frac1n- s \sum_{\substack{N<n\le KN\\ 
P(n)>\sqrt{s KN}}}\frac1n\\[5pt]
&= s \log K+O(s/N)- s\sum_{\sqrt{sKN}<p\le KN}\frac1p\sum_{N/p<m\le 
KN/p}\frac1m.
\end{align*}
Now,
\[
\sum_{N/p<m\le KN/p}\frac1m
=\begin{cases}
\log K+O(p/N),&p\le N\\
\log(KN/p)+O(1),&N<p\le KN.
\end{cases}
\]
Thus,
\begin{align*}
\sum_{\sqrt{sKN}<p\le KN} {\hskip -20pt}&{\hskip 
20pt}\frac1p\sum_{N/p<m\le KN/p}\frac1m\\[5 pt]
&=\sum_{\sqrt{sKN}<p\le N}\(\frac{\log K}{p}+O(1/N)\)+
\sum_{N<p\le KN}\(\frac{\log K}{p}+\frac{\log N-\log p+O(1)}{p}\)\\[5 pt]
&=\sum_{\sqrt{sKN}<p\le KN}\frac{\log K}{p}+\sum_{N<p\le 
KN}\frac{\log N-\log p}{p}+O(\log K/\log N)\\[5 pt]
&=\log 2\log K + O\bigg(  \dfrac{\log K \log(sK)}{\log N}  \bigg)
=  (\log 2+o(1))\log K.
\end{align*}
Hence, since $-\log\alpha(C_0)\le \sum_{n\in S(C_0)}1/n+O(s/N)$, we have
\[
-\log \a(C_0) \le s \big(1-\log2+o(1)\big)\log K\le
\big(1-(1-\log2)\eps+o(1)\big)\log L(N,s).
\]
Let $Q=L(N,s)^{1-\lambda}$, where $\lambda=\frac14((1-\log 2)\eps-4b^2)$.
  Also let $C'=\{(n,r)\in C_0:P(n)\le Q\}$.  As
before, using Lemma \ref{smoothsum} yields
$$
\del(C') = 1 + O\Bigg( s \sum_{\substack{n > N \\ P(n) \le Q}} 
\dfrac{1}{n}  \Bigg)
= 1 +o(1) \qquad (N\to \infty).
$$
Hence,
$$
\a(C_0)^{(1+1/Q)/\del(C')} \gg L(N,s)^{-1+(1-\log 2)\eps - \lambda}.
$$
On the other hand,
$$
\frac{s^{2} \log^2(QK)}{Q} \ll L(N,s)^{-1+4b^2+2\lambda}.
$$
By Lemma \ref{combinedlemma},
we have $\delta(C_0)>0$ for $N$ sufficiently large.
Thus, $\delta(C)>0$ by Lemma \ref{lem2pt1}.
\end{proof}

We now show that if $K$ is a bit smaller than in Theorem \ref{thmpos},
then in fact
$$
\delta(C)\ge(1+o(1))\alpha(C).
$$
The following result generalizes
Theorem B from the introduction.

\begin{thm}\label{theoremB1}
Suppose $0<\eps<1/2$, $0<b<\frac12 \sqrt{\eps}$ and $N\ge 100$.
Suppose that $C$ is a residue system with $S(C)$ consisting of
integers from $(N,KN]$ with multiplicity at most $s$,
where $s \le \exp\big(b\sqrt{\log N \log \log N}  \big)$
and $K =  L(N,s)^{( 1/2-\eps)/s}$.
Then
$$
\delta(C) \ge \(1+O\pfrac{1}{(\log N)^{\lambda}} \) \alpha(C),
$$
where $\lambda$ is a
positive constant depending only on $\eps$ and $b$.
\end{thm}

\begin{proof}
We follow the same general plan as in the proof of Theorem \ref{thmsmallsum}.
Since the sum of $1/n$ for all $n\in(N,KN]$ is $\log K+O(1/N)$ we have
\[
\alpha(C)\gg L(N,s)^{-1/2+\eps}.
\]
Let $Q=L(N,s)^{1/2-\lambda}$, where $\lambda = \frac13(\eps-4b^2)$.
In particular $Q\ge \log^2 N$.
Let $u=\log N/\log Q$, and let $C'$ be as in
Lemma \ref{combinedlemma}.  By Lemma \ref{smoothsum}, we have
\[
1-\delta(C')\ll\frac{s\log Q}{u^u}\ll\frac{s\log N}{(s\log N)^{2+\lambda}},
\]
so that $1/\delta(C')=1+O\((s\log N)^{-1-\lambda}\)$.
Since $|\log\a(C)|\le\log N$, we have
\[
\a(C)^{(1+1/Q)/\delta(C')}=(1+O(1/(\log N)^{\lambda} )\a(C).
\]
So, by Lemma \ref{combinedlemma} it suffices to show that
$s^2 (\log QK)^2/Q=O(\a(C)(\log N)^{-\lambda})$.  But, for large $N$ we have
$s^2 \le L(N,s)^{4b^2+\lambda}$.  Thus,
\[
\frac{s^2 (\log QK)^2}{Q}\ll\frac{s^2 \log^2 L(N,s)}{L(N,s)^{1/2-\lambda}}
\ll \frac{1}{L(N,s)^{1/2-2\lambda-4b^2}} \ll \frac{1}{L(N,s)^{1/2-\eps+\lambda}}
\ll \frac{\a(C)}{L(N,s)^{\lambda}}.
\]
This completes the proof.
\end{proof}

%
\section{Coverings and near-coverings of the integers}\label{near}
%

In this section, we address two items.  The first
shows that there are coverings of the integers with the moduli 
bounded below by $N$
and the multiplicity of the moduli \textit{near} the upper bound
on the multiplicity of the moduli given by Theorem \ref{thmsmallsum}.
The second shows that, when we allow $K$ to be large,
the density of the integers which are not covered
by a covering system using distinct moduli from $(N,KN]$ can be 
considerably smaller
than what is suggested by Theorem \ref{theoremB1}.

\begin{thm}\label{coveringexample}
For sufficiently large $N$ and $s=\exp(\sqrt{\log N\log \log N})$,
there exists an exact covering system with squarefree moduli greater 
than $N$ such that the multiplicity of each modulus does not exceed $s$.
\end{thm}

\begin{proof}
Let $p$ denote a prime and let $X_j=(j+1)^{j+1}$ for $j=0,1,\dots$.
We first show that
\begin{equation}
\label{Xineq}
\sum_{X_{j-1}<p\le X_j}[X_j/p]\ge X_{j-1}\quad (j\ge1).
\end{equation}
Here $[x]$ denotes the largest integer which is $\le x$.
Note that (\ref{Xineq}) holds for $j\le 5$.  Suppose then that
$j\ge6$.  Using the estimates (3.4), (3.17), and (3.18)
in Rosser and Schoenfeld \cite{RS}, we have that
\begin{align*}
\sum_{X_{j-1}<p\le X_j}[X_j/p]
&\ge X_j\sum_{X_{j-1}<p\le X_j}1/p-\pi(X_j)\\
&\ge X_j\left(\log\frac{\log X_j}{\log X_{j-1}}-\frac{1}{\log^2X_{j-1}}
-\frac{1}{\log X_j-\frac32}\right).
\end{align*}
The expression in the parentheses is
\begin{align*}
&\log\frac{(j+1)\log(j+1)}{j\log 
j}-\frac{1}{j^2\log^2j}-\frac{1}{(j+1)\log(j+1)-\frac32}\\
&\qquad>\frac1{j+1}\left(\frac{j+1}{j}-\frac{j+1}{2j^2}-\frac{j+1}
{j^2\log^2j}-\frac{1}{\log(j+1)-3/(2j+2)}\right)
>\frac{0.43}{j+1}.
\end{align*}
And $X_j=(j+1)j^j(1+1/j)^j>2.5(j+1)j^j$.  Thus,
\[
\sum_{X_{j-1}<p\le X_j}[X_j/p]
>\frac{2.5(j+1)0.43}{j+1}j^j>j^j,
\]
which proves (\ref{Xineq}).

We describe now an explicit construction of a covering
system, which we will then show satisfies the conditions of the theorem.
For $J\ge1$ and $s=X_J$, we establish that
there exists an exact covering system $C_{J}$ with squarefree moduli 
greater than
$$N_J=\prod_{j=0}^{J-1}X_j$$
such that the multiplicity of each modulus does not exceed $s$.
Set
$$P_j=\{p:\,X_{j-1}<p\le X_j\}.$$
We construct $C_{J}$, through induction on $J$, by choosing
moduli of the form $p_1\cdots p_J$ where each $p_j\in P_j$.
Observe that such a product $p_1\cdots p_J$ is necessarily $>N_{J}$.
One checks that $C_{1} = \{ (2,0), (2,1) \}$ satisfies the conditions 
for $C_{J}$ with $J = 1$.
Now, suppose that we have $C_{J}$ as above for some $J \ge 1$.
Thus, we have an exact
covering system $C_{J}$ with moduli of the form $p_1\cdots p_J$ where 
each $p_j\in P_j$.
Fix such a modulus $n=p_1\cdots p_J$, and let 
$(n,r_1),\dots,(n,r_t)$, with $t\le X_J$,
be the pairs of the form $(n,r)$ in $C_{J}$.  Let $q_1<q_2<\cdots$ be 
the complete list of primes
from $P_{J+1}$. To construct $C_{J+1}$, we replace each pair
$(n,r_i)$, $i\le[X_{J+1}/q_1]$, with the $q_{1}$ pairs
$(nq_1,r_i+n\mu)$, where $\mu = 0,\dots,q_1-1$.  Notice that the 
multiplicity of the modulus
$nq_1$ is at most $[X_{J+1}/q_1] q_1\le X_{J+1}$.  Next, we replace each pair
$(n,r_i)$, $[X_{J+1}/q_1]<i\le[X_{J+1}/q_1]+[X_{J+1}/q_2]$, with the 
$q_{2}$ pairs
$(nq_2,r_i+n\mu)$, where $\mu = 0,\dots,q_2-1$.  We proceed with this 
construction until
all the pairs $(n,r_1),\dots,(n,r_t)$ are replaced with
new pairs.  As $t \le X_{J}$, this will happen at some point by
\eqref{Xineq}.
This completes the inductive construction of our exact covering 
systems $C_{J}$.

To complete the proof of the theorem, it suffices to show that
$\log N_J\log\log N_J\ge\log^2s$ for large $J$.  Now
\[
\log N_J=\sum_{j=1}^Jj\log j\ge\int_1^Jt\log t\,dt
>\frac12J^2\log J-\frac14J^2,
\]
so that
\[
\log\log N_J>2\log J+\log\log J-\log 2+\log(1-1/(2\log J))
>2\log J+\log\log J-1,
\]
for $J\ge 7$.  Thus,
\[
\log N_J\log\log N_J>
J^2\log^2J+\frac12J^2\log J\left(\log\log J-1.5-\frac{\log\log 
J}{2\log J}\right)
>J^2\log^2J+3J\log^2J,
\]
for $J\ge350$.  But $\log^2s=(J+1)^2\log^2(J+1)<J^2\log^2J+3J\log^2J$
in the same range.  This completes the proof of the theorem.
\end{proof}

\begin{remark}
A more elementary proof, that does not use the estimates from
\cite{RS}, is possible by defining the sequence $X_j$ inductively
as the minimal numbers for which (\ref{Xineq}) holds.
\end{remark}

Suppose $s = 1$ and $N$, $KN$ are integers in Theorem \ref{theoremB1}.
Then $S(C)$ consists of distinct integers chosen from $(N,KN]$ so that
\[
\alpha(C) \ge \prod_{j=N+1}^{KN} \bigg(  1 - \dfrac{1}{j}  \bigg) = 
\dfrac{1}{K}.
\]
Thus, Theorem \ref{theoremB1} implies a lower bound
of approximately $1/K$
for any $\delta(C)$ with $S(C)\subseteq (N,KN]$,
provided $K$ is not too large.  It is clear that
the expression $1/K$ is not far from the truth,
since the argument of the introduction gives a residue system $C$ with
$\delta(C)\le 1/K$.  However,
we might ask about the situation when $K$ is large
compared to $N$.  The following result shows that
$\delta(C)$ can in fact be considerably smaller
than $1/K$ when $K$ is much larger than $N$.

\begin{thm}\label{thm3}
Suppose $N$ and $K$ are integers with $N\ge 1$ and $K$ sufficiently large.
Then there is some
residue system $C$ consisting of distinct moduli from $(N,KN]$ such that
\[
\delta(C) \le \frac{1}{K} \exp \( - \frac{\log K}{3N}\).
\]
\end{thm}

Before giving a proof of the above theorem, we give a lemma that will also
play a role in the next section.  For a set $T$ of positive
integers, we let $\CC(T)$ be the set of
residue systems $C$ with $S(C)=T$ and where $(n,r)\in C$ implies
$1\le r\le n$.  Also, define
$$
W(T) = \# \CC(T) = \prod_{n \in T} n.
$$

\begin{lem}\label{expecteddelta}
Let $T$ be a set of positive integers.
Then the expected value of
$\delta(C)$ over $C\in\CC(T)$, denoted $\mathbf{E} \delta(C)$, is $\prod_{n\in T}(1-1/n)$.
\end{lem}

\begin{proof}  Put $W=W(T)$ and
say $1\le m\le W$.  The number of systems $C\in\CC(T)$
with $m\in R(C)$
is $\prod_{n\in T}(n-1)$, since for each $n\in T$, there are $n-1$ choices
for $r$ with $1\le r\le n$ and $r\not\equiv m\pmod n$.
Thus,
\[
\sum_{C\in\CC(T)}\delta(C)
=\sum_{C\in\CC(T)}\frac1W\sum_{\substack{1 \le m \le W\\[1pt] m\in R(C)}} 1
=\frac1W\sum_{m=1}^W\sum_{\substack{C\in\CC(T)\\ m\in R(C)}} 1
=\frac1W\sum_{m=1}^W\prod_{n\in T}(n-1)=\prod_{n\in T}(n-1).
\]
The result follows by dividing this equation by $W$.
\end{proof}

\begin{remark}
It is not hard to prove a version of Lemma~\ref{expecteddelta} that allows
for taking moduli from $T$ with multiplicity greater than~1.
\end{remark}

\begin{proof}[Proof of Theorem \ref{thm3}]
There is a covering system with distinct moduli and smallest modulus $25$
(a result of Gibson~\cite{G1}),
so Theorem \ref{thm3} follows for $N\le 24$.
Henceforth we may assume that $N\ge 25$; however our argument holds for $N\ge4$.
We shall construct a residue system $C=\{(n,r(n)):N<n\le KN\}$ as follows.
We will randomly choose the values of $r(n) \in [1,n]$ for $N < n \le 2N$
so that each residue class modulo $n$ is taken with the same probability $1/n$
and the variables $r(n)$ are independent.
Based on the random choice of such $r(n)$ for $N < n \le 2N$,
we then select the remaining values of $r(n)$ with
$2N < n \le KN$ via a greedy algorithm.  In fact,
we show that, under our construction, the expected value of $\delta(C)$
over all randomly chosen values of $r(n) \in [1,n]$ for $N < n \le 2N$ is
\[
  \le \frac{1}{K} \exp \( - \frac{\log K}{3N}\).
  \]
The result thus follows.

Let $C_{2N}=\{(n,r(n)):N<n\le 2N\}$, where each $r(n)$ is chosen randomly
from $[1,n]$.  From Lemma \ref{expecteddelta}, it follows that
$\mathbf{E} \delta(C_{2N})= 1/2$.  Hence,
by the arithmetic mean--geometric mean inequality,
\be\label{logdelta}
\mathbf{E} \log\delta(C_{2N}) \le -\log2.
\ee
We will also make use of Lemma \ref{expecteddelta} in another way.
If $D$ is a subset of the integers in $(N, 2N]$ and $\tC = \{ (d, 
r(d)): d \in D \}$,
then it is not difficult to see that the expected value of $\delta(\tC)$ over
all randomly chosen values of $r(d) \in [1,d]$ for $d \in D$ is the same as
the expected value of $\delta(\tC)$ over
all randomly chosen values of $r(n) \in [1,n]$ for $n \in (N, 2N]$; 
in other words,
the random selection of extra residue classes not associated with 
$\tC$ will not
affect the expected value $\delta(\tC)$.  Thus, Lemma 
\ref{expecteddelta} implies
$\mathbf{E} \delta(\tC)= \alpha(\tC)$ where the expected value is 
over all randomly
chosen $r(n) \in [1,n]$ for $N < n \le 2N$.

Suppose then that the values of $r(n) \in [1,n]$ for $N < n \le 2N$ have been
chosen randomly.
For $2N < j \le KN$, we describe how to select $r(j)$.
For this purpose, we set $C_j= \{ (n,r(n)):N<n\le j\}$.
We use the greedy algorithm to choose $r(j)$
to be a residue class modulo $j$ containing the largest proportion of
$R(C_{j-1})$.  As in the introduction, this gives trivially
\[
\delta(C_j) \le \(1-\frac{1}{j}\) \delta(C_{j-1}).
\]
We can sometimes do better.  If $j$ has a divisor $d$ with $N<d\le2N$,
then there are residue classes modulo $j$ not intersecting
$R(C_{j-1})$.  In particular, the residue class $r(d)~({\rm mod}~d)$ contains
$r~({\rm mod}~j)$ when $r\equiv r(d)\pmod d$.
Let
\[
D(j)=\{d:\, d|j ,N<d\le2N\},\quad \tC_j=\{(d,r(d)):d\in D(j)\}.
\]
Let $f(j)$ be the number of residue classes $r~({\rm mod}~j)$ for which
$r\not\equiv r(d)\pmod d$ for each $d\in D(j)$.  If we choose $r(j)$
appropriately from among these $f(j)$ choices for $r$, we have
\be\label{better}
\delta(C_j) \le \(1-\frac{1}{f(j)}\) \delta(C_{j-1}).
\ee
The last equality is nonsense if $f(j)=0$, but in that case we have
$R(C_{j-1})=\emptyset$, and the theorem is trivial.  Also, there is
nothing to prove
if $f(j)=1$ since then $R(C_j)=\emptyset$.  Throughout the following
we assume that
$f(j)>1$.

We see from \eqref{better} and linearity of expectation that
\be\label{Efj}
\mathbf{E}\log\delta(C_j)- \mathbf{E}\log\delta(C_{j-1})
\le \mathbf{E}\log\left(1-\frac{1}{f(j)}\right)
\le-\mathbf{E}\left(\frac{1}{f(j)}\right).
\ee
Using Lemma \ref{expecteddelta} as described above, we have
\[
\mathbf{E}\delta\big(\tC_j\big)=\prod_{d\in D(j)}\(1-\dfrac{1}{d}\).
\]
Since $j$ is a common multiple of the members of $D(j)$,
it follows that $\delta\big(\tC_j\big)=f(j)/j$, so that
\[\mathbf{E} f(j)=j\mathbf{E}\delta\big(\tC_j\big)=j\prod_{d\in 
D(j)}\(1-\dfrac{1}{d}\).\]
By the arithmetic mean-harmonic mean inequality, we thus have
\[
\mathbf{E}\left(\frac{1}{f(j)}\right) \ge j^{-1}\prod_{d\in 
D(j)}\(1-\dfrac{1}{d}\)^{-1}
\ge\frac1j + \sum_{d\in D(j)}\frac1{dj}.
\]
After substituting the last inequality into \eqref{Efj}, we get
\[
\mathbf{E}\log\delta(C_j)- \mathbf{E}\log\delta(C_{j-1})
\le-\frac1j - \sum_{d\in D(j)}\frac1{dj}.
\]
Thus,
\begin{align*}
\mathbf{E}\log\delta(C)-\mathbf{E}\log\delta(C_{2N})
&\le - \sum_{j=2N+1}^{KN}\frac1j - \sum_{j=2N+1}^{KN}\sum_{d\in
D(j)}\frac1{dj}\\[5 pt]
&= - \sum_{j=2N+1}^{KN}\frac1j - \sum_{d=N+1}^{2N}\sum_{2N/d<l\le KN/d}
\frac1{d^2l} \\[5 pt]
&= -\log(K/2) + O(1/N) - \sum_{d=N+1}^{2N}\frac{\log K+O(1)}
{d^2}
\end{align*}
We have for $N\ge 4$ the estimate
$$
\sum_{d=N+1}^{2N}\frac1{d^2} \ge \int_{N+1}^{2N+1}
\frac{dt}{t^2} = \frac{N}{(N+1)(2N+1)} \ge \frac{1}{2.9N}.
$$
Therefore, by \eqref{logdelta},
$$
\mathbf{E}\log\delta(C) \le -\log K - \frac{\log K + O(1)}{2.9N}
$$
The theorem now follows.
\end{proof}

%
\section{Normal value of $\delta(C)$}\label{normal}
%

It is reasonable to expect that
$\delta(C)\approx \alpha(C)$ for almost all residue systems $C$ with fixed
$S(C)$.
In this section, we establish such a result when $S(C)$ consists of
distinct integers, by
considering the variance of $\delta(C)$ over $C\in\CC(T)$, where, as before,
$\CC(T)$ is the set of residue systems $C$ with $S(C)=T$.

\begin{thm}
Let $T$ be a set of distinct positive integers with minimum element $N
\ge 3$.  Let $\a$ be the common value of $\a(C)$ for $C\in\CC(T)$.  Then,
$$
\frac{1}{W(T)} \sum_{C \in \CC(T)} |\del(C)-\a|^2 \ll  \frac{\a^2\log N }{N^2}.
$$
\end{thm}

\begin{proof} From Lemma \ref{expecteddelta}, we have
${\bf{E}}\delta(C)=\a(C)=\a$.  Writing $W=W(T)$,
we deduce then that
\be\label{eq5pt1}
\frac{1}{W}\sum_{C\in\CC(T)}|\delta(C)-\a|^2=
\frac{1}{W}\sum_{C\in\CC(T)}\big(\delta(C)^2-\a^2\big).
\ee
We have
\[
\sum_{C\in\CC(T)}\delta(C)^2 
=\sum_{C\in\CC(T)}\bigg(\frac1W\sum_{\substack{1\le m\le W\\[1pt] 
m\in R(C)}}1\bigg)^2
=\frac1{W^2}\sum_{1\le m_1,m_2\le W}\sum_{\substack{C\in\CC(T)\\ 
m_1,m_2\in R(C)}}1.
\]
As in the proof of Lemma \ref{expecteddelta}, the inner sum is
\benn
\begin{split}
\prod_{\substack{n\in T\\ m_1\equiv m_2~({\rm mod}~n)}}{\hskip -15pt}&(n-1)
\prod_{\substack{n\in T\\ m_1\not\equiv m_2~({\rm mod}~n)}}{\hskip -15pt}(n-2)
=
W\prod_{n\in T}\(1-\frac2n\)\prod_{\substack{n\in T\\m_1\equiv 
m_2~({\rm mod}~n)}}
{\hskip -10pt}\frac{n-1}{n-2}\\
&=
\a^2W\prod_{n\in T}\frac{1-\frac2n}{\(1-\frac1n\)^2}
\prod_{\substack{n\in T\\ n\mid (m_1-m_2)}}\frac{n-1}{n-2}
=
\a^2W\prod_{n\in T}\(1-\frac1{(n-1)^2}\)
\prod_{\substack{n\in T\\ n\mid (m_1-m_2)}}\frac{n-1}{n-2}.
\end{split}
\eenn
Let $u=\sum_{n\in T}1/n^2$ and define
$ f(m_1,m_2)=\prod_{{n\in T,\  n\mid (m_1-m_2)}} (n-1)/(n-2)$.
Thus,
\be\label{keyident}
\sum_{C\in\CC(T)}\delta(C)^2
=
\frac{\a^2}{W}\(1-u+O\(\frac1{N^2}\)\)
\sum_{1\le m_1,m_2\le W}f(m_1,m_2).
\ee

For $S$ a finite set of integers which are $\ge 3$,
let $M(S)$ denote $\prod_{n\in S}(n-2)$, and let
$L(S)$ denote the least common multiple of the members of $S$.
We have
\[
f(m_1,m_2)=\prod_{\substack{n\in T\\ n\mid (m_1-m_2)}}\(1+\frac1{n-2}\)
=\sum_{\substack{S\subseteq T\\ L(S)\mid (m_1-m_2)}}\frac1{M(S)}.
\]
Thus,
\be\label{fsum}
\sum_{1\le m_1,m_2\le W}f(m_1,m_2)=\sum_{S\subseteq T}\frac1{M(S)}
\sum_{\substack{1\le m_1,m_2\le W\\ L(S)\mid (m_1-m_2)}}1
=W^2\sum_{S\subseteq T}\frac1{M(S)L(S)}.
\ee
In this last sum we separately consider the terms with $\# S\le 1$ and $\# S\ge2$.
We have
\be\label{importantterms}
\sum_{\substack{S\subseteq T\\ \# S\le 1}}\frac1{M(S)L(S)}
=1+\sum_{n\in T}\frac1{(n-2)n}
=1+u+O\(1/N^2\).
\ee
If $S\subseteq T$ and $\# S\ge2$, let $k>h$ be the largest two members of $S$.
Then $L(S)\ge\lcm[h,k]=hk/\gcd(h,k)$, so that
\[
E:=\sum_{\substack{S\subseteq T\\ \# S\ge2}}\frac1{M(S)L(S)}
\le\sum_{k>h\ge N}\frac{\gcd(h,k)}{(h-2)(k-2)hk}\sum_{U\subseteq [N,h-1]}
\frac1{M(U)}.
\]
The inner sum here is identical to $\prod_{N\le n\le h-1}(n-1)/(n-2)=(h-2)/(N-2)$,
so that
\[
E\ll\frac1N\sum_{k>h\ge N}\frac{\gcd(h,k)}{hk^2}
\le \frac1N\sum_{d\ge1}\sum_{\substack{k>h\ge N\\ d\mid h,~d\mid 
k}}\frac{d}{hk^2}
= \frac1N\sum_{d\ge 1}\sum_{v>w\ge N/d}\frac1{d^2wv^2}
\ll \frac1N\sum_{d\ge 1}\sum_{w\ge N/d}\frac1{d^2w^2}.
\]
In this last double sum,
if $d\le N$, then the sum on $w$ is $\ll d/N$, so that the contribution to $E$
is $\ll(\log N)/N^2$.  And if $d>N$, the sum on $w$ is $\ll 1$, so 
that the contribution
to $E$ is $\ll 1/N^2$.  We conclude that  $E\ll(\log N)/N^2$.
Thus, with \eqref{fsum} and \eqref{importantterms} we have
\[
\sum_{1\le m_1,m_2\le W}f(m_1,m_2)=W^2\(1+u+O((\log N)/N^2)\),
\]
so that from \eqref{keyident} and $u\ll 1/N$, we get
\[
\sum_{C\in\CC(T)}\delta(C)^2=\a^2W\(1+O((\log N)/N^2)\).
\]
The result now follows immediately from \eqref{eq5pt1}.
\end{proof}

\end{document}